\documentclass[a4paper,12pt]{article}
\usepackage{amsmath,amssymb,amsthm,geometry}
\usepackage{xcolor}
\usepackage{authblk}
\usepackage{setspace}
\usepackage{booktabs}
\usepackage{longtable}
\usepackage{array}
\usepackage[round,authoryear]{natbib} 
\usepackage{algorithm}
\usepackage{algpseudocode}
\usepackage{amsmath,amssymb,booktabs,longtable,array}
\usepackage{enumitem}
\usepackage[hidelinks]{hyperref}
\usepackage{graphicx}
\newcolumntype{L}[1]{>{\raggedright\arraybackslash}p{#1}}
\newcommand{\vect}[1]{\left[#1\right]}

\usepackage{caption}
\title{Preference-Based Estimation of Achievement Scalarising Function Parameters via Ordinal Regression}
\author[1]{Maria Barbati\thanks{Corresponding author. Email: \texttt{maria.barbati@unive.it}}}
\author[2]{Salvatore Greco}
\affil[1]{Department of Economics, Ca' Foscari University of Venice, Venice, Italy}
\affil[2]{Department of Economics and Business, University of Catania, Catania, Italy}
\date{}

\begin{document}
\maketitle
\begin{abstract}
Estimating the parameters of an achievement scalarising function is a challenging task in multiobjective optimisation. Achievement scalarising functions are widely used in several contexts, including preference-based optimisation methods and algorithmic procedures for generating or selecting efficient solutions. However, their effectiveness strongly depends on the appropriate specification of their parameters, such as weights and reference points.

In this paper, we propose an ordinal-regression-based methodology to estimate the parameters of an achievement scalarising function from preference information provided by the DM. In particular, the proposed approach allows us to infer both the weights associated with the objectives and the reference points characterising the scalarising function. The methodology is illustrated through a didactic example based on a multiobjective knapsack problem, which shows step by step how the preference information is translated into model parameters. 
The example illustrates that ordinal regression provides a flexible and interpretable framework for calibrating achievement scalarising functions and supporting the elicitation of decision-maker preferences in multiobjective optimisation.

\end{abstract}
\section{Introduction}

A classical way to incorporate preference information in multiobjective
optimisation is to transform the multiobjective function into a scalar one by
means of a scalarising function. Among the most important scalarising approaches,
achievement scalarising functions play a central role 
\citep{wierzbicki1980use}. The idea is to represent the
aspirations of the DM through a reference point in the objective space and to measure the achievement of each feasible solution with respect to such a point.
This makes it possible to search for solutions that are close to, or improve
upon, the desired aspiration levels. Achievement scalarising functions (ASF) are particularly attractive because, under suitable assumptions, they can generate Pareto optimal solutions, including unsupported efficient solutions that cannot be obtained by a simple weighted-sum scalarisation \citep{miettinen2002scalarising,nikulin2012new}.

A crucial issue in the use of such a function is the elicitation of its
parameters, namely the weights and the reference point. In many applications, asking the DM to provide these parameters directly may generate a significant cognitive burden \citep{klamroth2008integrating}.
Thus, there emerges the opportunity of estimating not only the
weights, but also the aspiration levels defining the reference point. Indeed, reference points are usually interpreted as desirable values for the objective functions and are often provided or progressively revised by the DM in interactive reference point methods \citep{luque2009incorporating}. Moreover,
the idea that a single fixed reference point may be too restrictive has already appeared in the literature. For example, multiple reference points have been used to approximate different regions of the Pareto front
\citep{figueira2010parallel}, while recent developments on bivariate achievement scalarising functions allow the reference point itself to be treated as a variable chosen within a set of admissible reference points
\citep{DECASTRO2026107469}. To the best of our knowledge, however, the
possibility of inferring such aspiration levels from holistic preference
information has not been previously investigated. 

For this reason, we adopt an indirect preference elicitation approach based on ordinal regression. Ordinal regression was introduced in multiple criteria decision aiding through the UTA method proposed by Jacquet-Lagrèze and Siskos \citep{jacquet1982assessing} to infer parameters such as weights from holistic preference information. The ordinal regression paradigm has been further developed in several directions, such as the Robust Ordinal
Regression, which considers the whole set of preference model parameters
compatible with the preference information provided by the DM, rather than a
single representative model \citep{greco2008ordinal}. The main principles and applications of Robust Ordinal Regression in multiple criteria decision
aiding are discussed in \citep{greco2010robust}, while further developments have extended the approach to preference learning and ranking problems
\citep{corrente2013robust}, as well as to value functions handling interacting criteria \citep{greco2014robust}.
More recently, ordinal regression has also been combined with richer preference elicitation procedures, i.e., the deck of cards method, \citep{barbati2024deck}
and represents one of the most popular methods in multiple criteria decision analysis \citep{greco2025fifty}.

The main contribution of this paper is to use ordinal regression to infer both the weights and the reference point of an ASF. More precisely, given a multiobjective optimisation problem, we use preference statements provided by the DM on a set of alternatives to estimate the parameters of the ASF. In this way, the ASF is not fixed a priori but is learned from the preferences expressed by the DM. 

The proposed approach has two main advantages. First, it reduces the cognitive effort required from the DM, who is asked to provide holistic preference information rather than precise numerical values of weights and aspiration levels. Second, it yields a scalarising model that can be directly embedded into a multiobjective optimisation procedure, thus providing a bridge between preference learning and optimisation. The resulting model can be used to guide the search toward regions of the Pareto front that are consistent with the Decision Maker (DM) preferences.

\textcolor{black}{Although the proposed approach is illustrated through a multiobjective knapsack problem, the framework is not specific to this application. The ordinal regression model is formulated directly in terms of a set of alternatives and their evaluations with respect to multiple criteria and can therefore be applied to general multiobjective optimisation problems. In particular, the criteria $g_j$, $j\in J$, can be interpreted either as performance criteria used to evaluate a finite set of alternatives or, more generally, as objective functions of a multiobjective optimisation problem. The role of the knapsack problem in this paper is therefore to provide a concrete and intuitive application in which the proposed approach can be illustrated and numerically tested.
}

The remainder of the paper is organised as follows. In Section ~\ref {sec:literature_review} the literature review is presented. Next, in Section ~\ref{sec:method} the ordinal regression model is formulated
 to infer the parameters of an achievement scalarising function. Section~\ref{sec:example} illustrates the proposed
procedure by adapting the approach for a multiobjective knapsack problem, proposing also a didactic example.  Finally, Section~\ref{sec:conclusions} provides some
concluding remarks and outlines possible directions for future research.

\section{Literature review} \label{sec:literature_review}

The theoretical foundations of scalarising functions were developed by \citep{wierzbicki1980use}. In this framework, the aspirations of the DM are represented by a reference point in the objective space, and feasible solutions are evaluated according to their achievement with respect to that point. More precisely, a main component
focuses on the objective that is farthest from the desired level, while the
augmentation term refines the comparison among solutions with the same worst
performance. The satisficing interpretation of this approach was further developed by \cite{wierzbicki1982mathematical}, where aspiration levels are interpreted as satisfactory target values rather than necessarily attainable optima. Later, \cite{wierzbicki1986completeness} studied the completeness and constructiveness of parametric characterisations of vector optimisation problems, providing further theoretical support for scalarising approaches.

The implementation of the reference-point approach had also been discussed by \cite{kallio1980implementation}, where the DM interactively specifies reference values which are then used to guide the search toward preferred efficient solutions. A comprehensive treatment of nonlinear multiobjective optimisation, scalarising functions, ideal and nadir points, and interactive methods is provided by \cite{miettinen1999nonlinear}. In addition, \cite{miettinen2002scalarising} analyse different scalarising functions and emphasise that ASFs can generate Pareto-optimal solutions even in non-convex cases, unlike weighted-sum functions.

Several authors have proposed variants or extensions of ASFs. \cite{ruiz2008additive} introduce an additive ASF for multiobjective programming problems that minimises the sum of the positive weighted deviations from the reference values, while \cite{nikulin2012new} propose a new ASF  introducing an integer parameter that controls the aggregation of positive weighted deviations from the reference
point. \cite{luque2012two} develop a two-slope ASF distinguishing between achievable and unachievable reference points.

ASFs have also been used in algorithmic frameworks, such as evolutionary optimisation algorithms \citep{ruiz2015preference}. Moreover, \cite{ishibuchi2010indicator} use ASFs to approximate the hypervolume indicator in evolutionary many-objective optimisation.

A relevant issue when adopting a scalarising function is the definition of the reference point, which is usually interpreted as a vector of desirable objective values. While the ideal and nadir points may help the DM understand the range of attainable objective values, the reference point expresses the DM's aspiration levels \citep{miettinen1999nonlinear, miettinen2002scalarising}. However, the idea that a single reference point may be too restrictive appears in several works. \cite{figueira2010parallel} propose a parallel multiple-reference-point approach for multiobjective optimisation, in which several reference points are used to explore different regions of the Pareto front. In this case, reference points are not estimated from the preferences of the DM but are used to generate a diverse approximation of the front. Similarly, \cite{aliano2021exact} propose a scalarisation method with multiple reference points for bi-objective integer linear optimisation problems. Their method uses multiple reference points algorithmically to generate Pareto-efficient solutions and compares this approach with modified Tchebycheff scalarisations.

A particularly important contribution is that of \cite{luque2015equivalent}, who introduced the concept of equivalent reference points in multiobjective programming. They show that, given a reference point, a vector of weights, and a nondominated solution obtained through an ASF, there exists a family of alternative reference points that generate the same solution.

Recent work also considers reference points as variable objects. \cite{DECASTRO2026107469} develop bivariate ASFs in which the reference point can be chosen within an admissible set. However, \cite{DECASTRO2026107469} do not infer the reference point from decision-maker preferences. Instead, the reference point becomes a variable of the scalarising framework. \cite{tanabe2024quality} review quality indicators for preference-based evolutionary multiobjective optimisation using a reference point and show that the position of the reference point has a strong impact on the evaluation of algorithms and on the induced region of interest.

ASFs are also extensively adopted in interactive multiobjective optimisation \citep{miettinen2008introduction}. In these methods, the DM participates in the optimisation process by providing aspiration levels, pairwise comparisons, or other preference statements. \cite{miettinen2006experiments} study classification-based scalarising functions, where the DM classifies objectives according to desired changes rather than specifying a single reference point. \cite{wierzbicki2007reference} further discusses reference-point approaches and objective ranking, showing how aspiration levels and objective rankings can be used to support interactive decision-making.

Among reference-point-based interactive methods, \cite{luque2009incorporating} incorporate preference information into interactive reference-point methods by allowing the DM to provide additional preference information beyond the reference point. Their approach recognises that a reference point alone may not fully capture the DM's preferences and that additional preference information can improve the projection of the reference point onto the Pareto front.
The same idea has been developed in the context of multiobjective integer programming by \cite{lokman2016interactive}, who propose an interactive algorithm to find the most preferred solution under an unknown quasiconcave value function. \cite{karakaya2018interactive} extend this line of research by developing interactive algorithms for a broad family of underlying preference functions. These approaches use pairwise preference information to reduce the set of potentially optimal solutions and to guide the search.

A particularly close line of research is based on weighted Tchebycheff preference functions. \cite{karakaya2021evaluating} study the evaluation of solutions and solution sets under multiple objectives using preference functions, including weighted Tchebycheff functions. \cite{karakaya2023finding} propose an interactive exact algorithm for multiobjective integer programmes under the assumption that the DM  preferences are consistent with a weighted Tchebycheff function. Their method presents pairs of solutions to the DM, updates the set of compatible weights, and uses the resulting information to search for the most preferred solution. However, they focus on estimating or restricting weights, whereas the present approach estimates both the weights and the reference point of an ASF. The interactive Tchebycheff framework has recently been extended by \citep{karakaya2025improved}, who address issues such as equal preferences and generalise the distance measure through a parameter that can represent different distance functions, including Euclidean and Chebyshev-type distances. 

The literature reviewed above shows that ASFs are well-established scalarising tools, that reference points may be multiple or non-unique, and that preference information can be used to guide the exploration of the Pareto front. However, most ASF-based approaches assume that the reference point is provided directly by the DM. Multiple-reference-point methods such as \citep{figueira2010parallel} and \citep{aliano2021exact} generate reference points algorithmically, while equivalent-reference-point results such as \citep{luque2015equivalent} analyse the non-uniqueness of reference points after a solution has been obtained. Interactive methods such as \citep{luque2009incorporating}, \citep{miettinen2006synchronous}, \citep{phelps2003interactive}, \citep{lokman2016interactive}, and \citep{karakaya2023finding} use decision-maker preferences to update the search, but they do not estimate the reference point of an ASF from holistic preference information.

\section{Ordinal regression for achievement scalarising functions}\label{sec:method}

A scalarising function assigns to each alternative $a \in A$ a value $U(a)$, based on the evaluations $g_j(a)$ of $a$ with respect to the considered criteria $g_j$, $j \in J$. Without loss of generality, the criteria are assumed to be represented by real-valued functions $g_j: A \rightarrow \mathbb{R}$, such that the greater $g_j(a)$ is, the better the evaluation of $a$ with respect to criterion $g_j$. 
 More precisely, each $a \in A$ is assigned the overall evaluation
\[
U(a)=
\rho \sum_{j\in J} w_j\left(g_j(a)-r_j\right)
-
(1-\rho)\max_{j\in J} w_j\left(r_j- g_j(a)\right)=
\]
\begin{equation}\label{value_function}
(1-\rho)\min_{j\in J} w_j\left(g_j(a)-r_j\right)
+
\rho \sum_{j\in J} w_j\left(g_j(a)-r_j\right),
\end{equation}
 where 
 \begin{itemize}
     \item $w_j, j \in J$, is the weights assigned to criterion $g_j, j \in J$,
     \item $r_j, j\in J$ is the evaluation of a reference point $\mathbf{r}=[r_1,\ldots,r_{|J|}]$ with respect to criterion $g_j$
     \item $0 \leq\rho \leq 1$ is a a fixed parameter, generally chosen close to zero.
     \end{itemize}

 With respect to the weights $w_j, j\in J$ it is assumed that $w_j\ge 0$ for all $j \in J$ and $\displaystyle{\sum_{j \in J} w_j=1}$. 
Let us note that the same framework can also be used for multiobjective problems by interpreting the criteria $g_j$ as objective functions. Thus, each feasible solution $a \in A$ is associated with a vector of objective values $g_j(a)_{j \in J}$, which can be transformed into a single scalar by the scalarising function $U$. \textcolor{black}{The framework described above is general and is not tied to a specific optimisation problem. In a multiobjective optimisation setting, the criteria $g_j$, $j\in J$, can be directly interpreted as objective functions, with each feasible solution $a\in A$ associated with the objective vector
\[
\mathbf{g}(a)=\left(g_1(a),\ldots,g_{|J|}(a)\right).
\]
The scalarising function $U$ then maps this objective vector into a single overall value, thereby inducing a preference ordering over the feasible solutions. Consequently, the proposed ordinal regression approach can be used to induce the parameters of an ASF from preference information provided by the DM for a broad class of multiobjective optimisation problems. The specific optimisation problem considered in the application only determines the feasible set $A$ and the objective functions $g_j$; the proposed preference elicitation and parameter induction framework remains unchanged.
}

\textcolor{black}{The following two remarks are useful for the ordinal regression model developed in the remainder of the paper.
First, the reference point affects the additive component of the ASF only through an alternative-independent constant. Indeed, defining
\begin{equation}\label{value_function_bis}
U^\ast(a)=
(1-\rho)\min_{j\in J} w_j\left(g_j(a)-r_j\right)+\rho \sum_{j\in J} w_j g_j(a)
\end{equation}
we have
\[
U(a)=U^\ast(a)-\rho\sum_{j\in J}w_jr_j.
\]
Since the second term is constant with respect to $a$, it follows that, for all $a',a''\in A$,
\[
U(a')\geq U(a'')
\quad\Longleftrightarrow\quad
U^\ast(a')\geq U^\ast(a'').
\]
Therefore, $U$ and $U^\ast$ induce exactly the same preference ordering over $A$. This observation implies that, when ordinal preferences are used to estimate the parameters of the ASF, the reference point affects the induced ordering only through the Tchebycheff component.
Second, the reference point inducing a given preference ordering is generally not unique. In particular, let $\mathbf{r}^\circ$ be another reference point such that, for all $j\in J$,
\[
w_jr_j^\circ=w_jr_j+k,
\]
for some $k\in\mathbb{R}$. Consider the corresponding ASF
\[
U^\circ(a)=
(1-\rho)\min_{j\in J}w_j\left(g_j(a)-r_j^\circ\right)
+
\rho\sum_{j\in J}w_j\left(g_j(a)-r_j^\circ\right).
\]
Then
\[
w_j\left(g_j(a)-r_j^\circ\right)
=
w_j\left(g_j(a)-r_j\right)-k,
\]
and hence
\[
U^\circ(a)
=
U(a)-\left[(1-\rho)+\rho|J|\right]k.
\]
The difference between $U^\circ(a)$ and $U(a)$ is therefore constant with respect to $a$. Consequently, for all $a',a''\in A$,
\[
U(a')\geq U(a'')
\quad\Longleftrightarrow\quad
U^\circ(a')\geq U^\circ(a'').
\]
Thus, $\mathbf{r}$ and $\mathbf{r}^\circ$ induce the same preference ordering over $A$. This result, established in \cite{luque2015equivalent}, highlights that the reference point identified from ordinal preference information should generally be interpreted as one representative of a family of reference points compatible with the observed preferences.
}

\subsection {Inducing the weights and the reference point of the ASF function}

In this paper we introduce an ordinal regression approach for ASF to induce the weights $w_j, j \in J$ and the reference point $\mathbf{r}=[r_1,\ldots,r_{|J|}]$ of a scalarising function from DM's preference information expressed in terms of binary comparisons ``$a$ is at least as good as $b$'', denoted by $a \succsim_{DM} b$, with $a,b \in A^R$ and $A^R$ is a set of reference alternatives. Observe that we put $v_j=w_j \cdot r_j$, so that $r_j=\frac{v_j}{w_j}, j\in J$.  
The weights $w_j, j\in J$ and the reference point $\mathbf{r}=[r_1,\ldots,r_{|J|}]$ that best represent the DM's preferences are obtained by solving the following MILP model, denoted \textbf{ORD-ASF}, with respect to the following unknown variables: 
\begin{itemize}
    \item $w_j, j\in J$,
    \item $v_j, j\in J$,
    \item $\sigma^+(a), \sigma^-(a), a\in A^R$,
    \item $t(a), a \in A^R$,
    \item $\gamma_j(a), j \in J, a \in A^R$:
\end{itemize}

\textbf{ORD-ASF}
\begin{align}
\min \quad 
& \sum_{a \in A^R} \left( \sigma^+(a) + \sigma^-(a) \right)
&&  \label{eq:ord-asf-obj} \\
\text{s.t.} \quad
& \sum_{j \in J} w_j = 1
&&   \label{eq:ord-asf-normalization}\\
& t(a) \ge v_j - w_j g_j(a),
&& \forall a \in A^R,\; j \in J
&&  \label{eq:ord-asf-upper-bound}\\
& v_j - w_j g_j(a) \ge t(a)-\gamma_j(a)M,
&& \forall a \in A^R,\; j \in J
&&  \label{eq:ord-asf-big-m}\\
& \textcolor{black}{U(a)=\rho \sum_{j \in J} w_j g_j(a) - (1-\rho)t(a)},
&& \forall a \in A^R
&&  \label{eq:ord-asf-value}\\
& U(a)+\sigma^-(a)-\sigma^+(a)
\ge
U(b)+\sigma^-(b)-\sigma^+(b),
&& \forall a,b \in A^R : a \succsim_{DM} b
&&  \label{eq:ord-asf-preference}\\
& \sum_{j \in J} \gamma_j(a) \le |J|-1,
&& \forall a \in A^R
&&  \label{eq:ord-asf-active-big-m}\\
& \sigma^+(a)\ge 0,\quad \sigma^-(a)\ge 0,
&& \forall a \in A^R
&&  \label{eq:ord-asf-slack-nonnegative}\\
& w_j \ge 0,
&& \forall j \in J
&& \label{eq:ord-asf-weight-nonnegative}\\
& {\color{black} w_j g_{j*}\le v_j \le w_j g_j^*},
&& \forall j \in J
&& \label{eq:ord-asf-v_j-non-nul}\\
& \gamma_j(a)\in \{0,1\},
&& \forall a \in A^R,\; j \in J.
&& \label{eq:ord-asf-binary}
\end{align}

The objective function \eqref{eq:ord-asf-obj} minimises the total deviation from the preference information provided by the DM. The variables $\sigma^+(a)$ and $\sigma^-(a)$  measure the possible inconsistencies between the preferences expressed by the DM and those represented by the model. 

Constraint \eqref{eq:ord-asf-normalization} imposes the normalisation of the criteria weights. 
The set of constraints \eqref{eq:ord-asf-upper-bound} introduces the auxiliary variable \textcolor{black}{$t(a)=\max_{j\in J} w_j\left(r_j-g_j(a)\right)$} as an upper bound on the weighted deviations from the reference point. For each alternative and each criterion, $t(a)$ must be at least as large as the corresponding weighted deviation. 
The set of constraints \eqref{eq:ord-asf-big-m} is used to linearise the maximum operator through a big-$M$ formulation,  where $M>0$ is a sufficiently large valid constant. The binary variable $\gamma_j(a)$ determines whether the corresponding deviation is allowed to be smaller than $t(a)$ or whether it must coincide with it.
The set of constraints \eqref{eq:ord-asf-value} defines the overall value $U(a)$ that combines a small weighted-sum component, controlled by the parameter $\rho$, with a term depending on the maximum weighted deviation represented by $t(a)$. \textcolor{black}{Observe that in fact $U(a)$ is defined here on the basis of  $U^\ast(a)$ in equation \ref{value_function_bis} of Section \ref{sec:method}, by observing that 
\[
U(a)=
\rho \sum_{j\in J} w_jg_j(a)
+
(1-\rho)\min_{j\in J} w_j\left(g_j(a)-r_j\right)=
\]
\[
\rho \sum_{j\in J} w_jg_j(a)
-
(1-\rho)\max_{j\in J} w_j\left(r_j-g_j(a)\right)=\rho \sum_{j\in J} w_jg_j(a)
- (1-\rho)t(a).
\]
} 
The set of constraints \eqref{eq:ord-asf-preference} translates the pairwise preference information provided by the DM into mathematical inequalities. Whenever the DM states that alternative $a$ is at least as good as alternative $b$, the model imposes that the adjusted value of $a$ is not smaller than the adjusted value of $b$.  Let us note that different types of preference statements can be represented by suitably
modifying constraint \eqref{eq:ord-asf-preference}. In fact, in the case of a weak preference \(a \succeq_{\mathrm{DM}} b\), the standard preference
constraint is sufficient. Instead, if the DM expresses a strict preference \(a \succ_{\mathrm{DM}} b\), a strictly positive
discrimination threshold \(\varepsilon>0\) is introduced, preventing two alternatives declared to be strictly
ordered by the DM from receiving the same adjusted utility value, and the following constraint is
imposed:
$
U(a)+\sigma^{-}(a)-\sigma^{+}(a)
\geq
U(b)+\sigma^{-}(b)-\sigma^{+}(b)+\varepsilon.
$
The set of constraints \eqref{eq:ord-asf-active-big-m} ensures that, for each reference alternative, at least one of the big-$M$ constraints is active. As a consequence, the auxiliary variable $t(a)$ is not only an upper bound, but it coincides with one of the weighted deviations, thus representing the maximum deviation.
Constraints \eqref{eq:ord-asf-slack-nonnegative} and \eqref{eq:ord-asf-weight-nonnegative} impose the non-negativity of the slack variables and of the weights, respectively. \textcolor{black}{Constraints \eqref{eq:ord-asf-v_j-non-nul} ensures that auxiliary variables $v_j$, introduced through the relation $v_j=w_jr_j$, satisfy the condition that  $v_j=0$ whenever $w_j=0$. Otherwise, the mathematical formulation could assign a non-zero value to $v_j$ even when the corresponding criterion has zero weight, which would make it impossible to recover a finite reference-point coordinate $r_j$ from $r_j=v_j/w_j$ and, more importantly, would allow the parameter $v_j$ to affect the scalarising function despite the corresponding criterion having zero Tchebycheff weight. In particular, $g_{j*}$ and $g_j^*$ represents two estimated lower and upper bounds, possibly the minimum and the maximum, for the objective function $g_j$.} Constraint \eqref{eq:ord-asf-binary} defines the binary nature of the variables $\gamma_j(a)$. Let us note that, once the parameters $v_j$, $j \in J$, are fixed, the term
$
\rho \sum_{j \in J} v_j
$
is constant with respect to the alternatives. Therefore, it does not affect either their
ranking or the optimal solution and can be omitted from the scalarising utility function without loss
of generality. In the didactic example presented later, we adopt this equivalent formulation. \textcolor{black}{It is also worth noting that the recovery of the reference-point component \(r_j\) from
$
r_j=\frac{v_j}{w_j}
$
requires \(w_j>0\). If \(w_j=0\), the value of \(r_j\) cannot be identified from the preference information, since the corresponding criterion does not contribute to the ASF through the product \(w_jr_j\).}

Let us highlight that preference information provided by the DM can be expressed also in cardinal form. In this case, each reference alternative $a\in A^R$ is associated with a numerical evaluation
$\nu(a)$, obtained for instance, by means of a deck-of-cards-based elicitation procedure \citep{barbati2024deck}. The scalarising function is then required to  reproduce these evaluations, up to a non-negative scaling factor $k$,
while allowing for positive and negative deviations. Then the constraint \ref{eq:ord-asf-preference} of the \textbf{ORD-ASF} model could be modified as
\[
U(a)-\sigma^+(a)+\sigma^-(a)=k\nu(a),
\qquad \forall a\in A^R.
\]
By integrating cardinal preference information, the proposed approach enhances
the model's ability to represent the DM's preferences while making the elicitation process more user-friendly and less cognitively demanding.
\subsubsection{\textcolor{black}{Refinement strategies of the \textbf{ORD-ASF model}}}
\textcolor{black}{The proposed ordinal regression framework can be extended in several directions. First, the parameter $\rho$, which is fixed in the current formulation, could also be estimated from the preference information provided by the DM. In this case, the objective would be to jointly determine the weights, the reference point, and the value of $\rho$ that minimise the inconsistency with the observed preferences. However, when $\rho$ is treated as a decision variable, the resulting optimisation problem is no longer a MILP, since products between $\rho$ and the other unknown parameters arise in the scalarising function.}

\textcolor{black}{A straightforward approach to address this issue is to consider a finite grid of candidate values for $\rho$. For each fixed value $\rho_k$ in the grid, the \textbf{ORD-ASF} model can be solved to determine the weights and the reference point that minimise the total deviation from the DM's preferences. The value of $\rho$ associated with the smallest optimal value of the objective function can then be selected as the best candidate. This procedure can subsequently be refined through an alternating optimisation scheme. Starting from the best value $\rho^{(0)}$ identified on the grid, the weights and the reference point are first estimated by solving the corresponding \textbf{ORD-ASF} model with $\rho^{(0)}$ fixed. Then, keeping the estimated weights and reference point fixed, $\rho$ is updated by minimising the total deviation from the DM's preferences over $\rho\in[0,1]$. The newly obtained value of $\rho$ is then fixed again, and the weights and the reference point are re-estimated by solving the \textbf{ORD-ASF} model. The procedure is repeated until no further improvement in the objective function is obtained. Since the joint problem is non-linear and, in general, non-convex, this alternating procedure should be regarded as a heuristic refinement procedure rather than as a method guaranteeing a globally optimal solution. Multiple initial values of $\rho$ or a sufficiently fine grid can be used to reduce the risk of convergence to a locally optimal solution.}

\textcolor{black}{A second possible extension concerns the discrimination threshold used to represent
strict preference statements. In the basic formulation, a positive value \(\epsilon>0\) can be fixed a priori and imposed as a minimum separation between the corrected utilities of
two alternatives such that \(a \succ_{DM} b\). Alternatively, \(\epsilon\) can be treated as a variable and maximised in a second optimisation stage. In this case, after
solving the \textbf{ORD-ASF} model and obtaining the optimal inconsistency level \(\sigma^*\), one may solve a
second optimisation problem whose objective is
to maximise $
 \epsilon
$ 
subject to the constraints $
U(a)+\sigma^{-}(a)-\sigma^{+}(a)
\geq
U(b)+\sigma^{-}(b)-\sigma^{+}(b)+\varepsilon.
\quad \forall a,b\in A^R:\; a\succ_{DM} b,
$
and
$
\sum_{a\in A^R}\left(\sigma^+(a)+\sigma^-(a)\right)\leq \sigma^*.
$
Thus, \(\epsilon\) represents the minimum utility distance between strictly
ordered alternatives, and maximising \(\epsilon\) means maximising the worst-case
separation compatible with the same level of preference inconsistency.}

\textcolor{black}{A third possible extension consists in allowing the Tchebycheff and weighted-sum components of the ASF to use different vectors of weights. In particular, taking into account formulation (\ref{value_function_bis}), the scalarising function could be written as
\begin{equation}\label{Ext}
U(a)=
(1-\rho)\min_{j\in J} w_j^{T}\left(g_j(a)-r_j\right)
+
\rho\sum_{j\in J}w_j^{S}g_j(a),
\end{equation}
where $w^T$ and $w^S$ denote, respectively, the weights associated with the Tchebycheff and weighted-sum components. This formulation allows the relative importance assigned by the DM to the criteria to differ according to whether the preference model is driven by the worst-performing criterion or by the overall compensatory evaluation of the criteria.}

\textcolor{black}{In this more general preference model, the parameter $\rho$, with $0 \leqslant \rho < 1$, does not necessarily need to be estimated separately. Indeed, $\rho$ can be absorbed into the scale of the two weight vectors. More precisely, taking $U^{ext}(a)=\frac{U(a)}{1-\rho}$, we get $U(a') \geqslant U(a'')$ if and only if $U^{ext}(a') \geqslant U^{ext}(a'')$. Setting 
\begin{equation}\label{multiplier}
\mu=\frac{\rho}{1-\rho},
\end{equation}
we obtain
\[
U^{ext}(a)=\frac{1}{1-\rho}U(a)=
\]
\begin{equation}\label{Ext1}
\frac{1}{1-\rho}\left[(1-\rho)\min_{j\in J} w_j^{T}\left(g_j(a)-r_j\right)
+
\rho\sum_{j\in J}w_j^{S}g_j(a)\right]=\min_{j\in J} w_j^{T}\left(g_j(a)-r_j\right)
+
\mu\sum_{j\in J}w_j^{S}g_j(a),
\end{equation}
and, defining
\[
\widetilde{w}_j^{S}=\mu w_j^{S}, j=1,\ldots,m
\]
the scalarising function can be equivalently written as
\[
U^{ext}(a)=
\sum_{j\in J}\widetilde{w}^{S}_j g_j(a) +\min_{j\in J} w^{T}_j\left(g_j(a)-r_j\right)
\]
\begin{equation}\label{Ext2}
=\sum_{j\in J}\widetilde{w}^{S}_j g_j(a)-\max_{j\in J} w^{T}_j\left(r_j-g_j(a)\right).
\end{equation}
Assuming $\sum_{j\in J}w^T_j=1$ and $\sum_{j\in J}w^S_j=1$, we have 
\begin{equation}\label{Ext3}
\sum_{j\in J}\widetilde{w}^{S}_j=\sum_{j\in J}\mu w^{S}_j=\mu\sum_{j\in J} w^{S}_j=\mu 
\end{equation}
so that, remembering (\ref{multiplier}), from which we obtain 
\[
\rho=\frac{\mu}{1+\mu},
\]
finally we obtain 
\begin{equation}\label{Ext4}
\rho=\frac{\sum_{j\in J}\widetilde{w}^{S}_j}{1+\sum_{j\in J}\widetilde{w}^{S}_j}.
\end{equation}
Consequently, adopting the formulation (\ref{Ext}) of the scalarising function in \textbf{ORD-ASF} one can replace equation
\[
U(a)=\rho \sum_{j \in J} w_j g_j(a) - (1-\rho)t(a)
\] 
with equation 
\[
U(a)=\sum_{j \in J} \widetilde{w}^S_j g_j(a) - \widetilde{t}(a)
\]
where $\widetilde{t}(a)=\max_{j\in J} w^{T}_j\left(r_j -g_j(a)\right)$. Moreover, constraint (\ref{eq:ord-asf-normalization}) can be replaced by 
\[
\displaystyle{\sum_{j=1}^m w_j^T=1}
\]
while constraint (\ref{eq:ord-asf-weight-nonnegative}) becomes
\[
\widetilde{w}^S_j \ge 0, w_j^T \ge 0, j \in J.
\]
The \textbf{ORD-ASF} problem so reformulated can be solved for the unknown variables $\widetilde{w}^S_j$ and $w_j^T, j \in J$, so that, using (\ref{Ext3}) and (\ref{Ext4}), also the value of $\rho$ can be obtained.
}
\textcolor{black}{Given the weighted reference-point parameters as
$
v_j^T = w_j^T r_j$ for each criterion $j \in J$, it is worth noting that the recovery of the reference point from
\[
r_j=\frac{v_j^T}{w_j^T}
\]
requires $w_j^T>0$. If $w_j^T=0$, the corresponding term $w_j^T(r_j-g_j(a))$ is identically equal to zero for every alternative $a$ and for any value of $r_j$. Hence, in this case, the value of $r_j$ cannot be identified from the preference information, as it has no effect on the scalarising function or on the induced preference ordering. This does not affect the identification of the scalarising function itself, since the corresponding criterion does not contribute to the Tchebycheff component when $w_j^T=0$.}

\textcolor{black}{Accordingly, the reference point should be interpreted as being identified only for those criteria receiving a strictly positive Tchebycheff weight. For criteria with $w_j^T=0$, the value of $r_j$ is immaterial and can be assigned arbitrarily without affecting the representation of the DM's preferences. In practical implementations, one may either report the reference point only for criteria with $w_j^T>0$ or impose a small positive lower bound on the Tchebycheff weights if a fully specified reference point is required for all criteria.}

\textcolor{black}{This extension provides a more flexible representation of the DM's preferences, as it allows the importance of the criteria to differ between the non-compensatory Tchebycheff component and the compensatory weighted-sum component. It also raises an interesting empirical question, namely whether the additional flexibility provided by two distinct weight vectors leads to a significantly better representation of the DM's preference information than the standard ASF formulation based on a common weight vector.}

\subsection{Balancing the Induced Weights}

It may happen that one or more of the induced weights are equal to zero or that the DM aspire to more balanced weights. 
In this case, a second MILP can be formulated in order to avoid, whenever possible, null weights. 
More precisely, let
\[
\sigma^\ast =
\min \sum_{a \in A^R} \left( \sigma^+(a) + \sigma^-(a) \right)
\]
be the optimal value of the objective function obtained by solving the \textbf{ORD-ASF} model. 
Then, among all solutions that preserve the same minimum level of preference inconsistency, we select the one that maximises the minimum weight. 
This can be done by solving the following MILP, denoted as \textbf{ORD-ASF-WB}:

\textbf{ORD-ASF-WB}
\begin{align}
\max \quad 
& \omega
&&  \label{eq:ord-asf-wb-obj} \\
\text{s.t.} \quad
& \sum_{j \in J} w_j = 1
&&   \label{eq:ord-asf-wb-normalization}\\
& w_j \ge \omega,
&& \forall j \in J
\label{eq:ord-asf-wb-weight-lower-bound}\\
& t(a) \ge v_j - w_j g_j(a),
&& \forall a \in A^R,\; j \in J
\label{eq:ord-asf-wb-upper-bound}\\
& v_j - w_j g_j(a) \ge t(a)-\gamma_j(a)M,
&& \forall a \in A^R,\; j \in J
\label{eq:ord-asf-wb-big-m}\\
& U(a)=\rho \sum_{j \in J} w_j g_j(a) - (1-\rho)t(a),
&& \forall a \in A^R
\label{eq:ord-asf-wb-value}\\
& U(a)+\sigma^-(a)-\sigma^+(a)
\ge
U(b)+\sigma^-(b)-\sigma^+(b),
&& \forall a,b \in A^R : a \succsim_{DM} b
\label{eq:ord-asf-wb-preference}\\
& \sum_{j \in J} \gamma_j(a) \le |J|-1,
&& \forall a \in A^R
\label{eq:ord-asf-wb-active-big-m}\\
& \sum_{a \in A^R} \left( \sigma^+(a) + \sigma^-(a) \right)=\sigma^*,
&&
\label{eq:ord-asf-wb-optimal-inconsistency}\\
& \sigma^+(a)\ge 0,\quad \sigma^-(a)\ge 0,
&& \forall a \in A^R
\label{eq:ord-asf-wb-slack-nonnegative}\\
& w_j \ge 0,
&& \forall j \in J
\label{eq:ord-asf-wb-weight-nonnegative}\\
& \gamma_j(a)\in \{0,1\},
&& \forall a \in A^R,\; j \in J,
\label{eq:ord-asf-wb-binary}\\
& \omega \ge 0.
&&
\label{eq:ord-asf-wb-omega-nonnegative}
\end{align}

The model \textbf{ORD-ASF-WB} preserves the structure of \textbf{ORD-ASF}, but introduces an additional objective aimed at avoiding null weights whenever possible. The meaning of the constraints already included in \textbf{ORD-ASF} remains unchanged.
The objective function \eqref{eq:ord-asf-wb-obj} of \textbf{ORD-ASF-WB} maximizes the auxiliary variable $\omega$, which represents a common lower bound for all criteria weights. Therefore, the model selects, among the solutions that are optimal with respect to the preference inconsistency, the one in which the smallest criterion weight is as large as possible.
Constraint \eqref{eq:ord-asf-wb-weight-lower-bound} imposes that every criterion weight must be at least equal to $\omega$. Since $\omega$ is maximized, this constraint pushes all weights away from zero whenever this is compatible with the preference information provided by the DM.
Constraint \eqref{eq:ord-asf-wb-optimal-inconsistency} ensures that the solution of \textbf{ORD-ASF-WB} preserves the optimal level of preference inconsistency obtained by solving \textbf{ORD-ASF}, where $\sigma^*$ denotes the optimal objective value of \textbf{ORD-ASF}. In other words, the second model does not worsen the representation of the DM's preferences; it only selects, among the optimal solutions of the first model, one with more balanced weights.
Finally, constraint \eqref{eq:ord-asf-wb-omega-nonnegative} imposes the non-negativity of the lower bound on the weights.

\subsection{Preference-Based Reference Point Induction}

In the literature, it has emerged that the reference point used in an
ASF is not necessarily unique \citep{luque2015equivalent}. Namely, given a reference
point, a weight vector and an efficient solution obtained by solving a multiobjective optimisation problem with an ASF formulation, they prove that there exists a set of equivalent reference points yielding the same efficient solution. Therefore, the reference point inferred from preference information should not necessarily be interpreted as the unique representation of the DM aspirations, but rather as one representative of a family of preference-compatible reference points. \textcolor{black}{An additional source of non-identification arises when a criterion receives zero weight. If $w_j=0$, then
\[
w_j\left(g_j(a)-r_j\right)=0
\]
for every alternative $a$ and for any finite value of $r_j$. Hence, the corresponding coordinate $r_j$ has no effect on the scalarising function and cannot be inferred from the DM's preference information. In this case, any value assigned to $r_j$ leads to the same preference representation. Consequently, the set of reference points compatible with the observed preferences may be substantially larger than the one-dimensional family generated by the equivalence transformation above: in addition to the common translation-type degree of freedom, the coordinates associated with zero-weight criteria are completely unrestricted.
Therefore, when all Tchebycheff weights are strictly positive, the equivalence result characterises a one-dimensional family of preference-equivalent reference points. When one or more weights are zero, additional degrees of freedom arise, and the corresponding reference-point coordinates are not identified. More generally, if $q$ criteria have zero Tchebycheff weight, the set of preference-equivalent reference points contains $q$ additional unconstrained directions. The reference point should therefore be interpreted as a preference-compatible representation, with its coordinates being meaningfully identified only for criteria with strictly positive Tchebycheff weights.}

 Following that, for the model \textbf{ORD-ASF} the additional constraint
\[
\sum_{j\in J} v_j = 1
\]
can be included to normalise the parameters $v_j$, $j\in J$. Since $v_j=w_jr_j$, this corresponds to imposing
\[
\sum_{j\in J} w_jr_j=1.
\]
Thus, the constraint fixes the scale of the weighted reference point components.

Furthermore,  the set of reference points induced can be transformed to obtain a representative of the set of reference points such that their meaning can be clear to the DM.
Let \(\widehat v_j\), \(j \in J\), denote the value of parameter \(v_j\)
obtained from the \textbf{ORD-ASF} or \textbf{ORD-ASF-WB} model, and let
\(\widehat w_j\), \(j \in J\), be the criterion weights estimated by the
same model. Moreover, let \(g^*_j\) denote the ideal value of criterion
\(j\). In order to calibrate the values \(\widehat v_j\) with respect to
the ideal-point-based quantities \(\widehat w_j g^*_j\), we introduce a
common shifting variable \(\kappa\) and two non-negative deviation
variables \(\delta^+_j\) and \(\delta^-_j\), for each \(j \in J\). The
variable \(\kappa\) allows all values \(\widehat v_j\) to be shifted by
the same amount, while \(\delta^+_j\) and \(\delta^-_j\) measure,
respectively, the positive and negative discrepancies between the
shifted value \(\widehat v_j+\kappa\) and the target value
\(\widehat w_j g^*_j\).

The following model, denoted \textbf{ORD-ASF-Vcal}, is considered:
\begin{align}
\min \quad
& \sum_{j\in J} \left(\delta_j^+ + \delta_j^-\right)
&& \label{eq:ord-asf-vcal-obj} \\
\text{s.t.} \quad
& \widehat v_j + \kappa - \delta_j^+ + \delta_j^-
=
\widehat w_j g_j^*
&& \forall j\in J
\label{eq:ord-asf-vcal-target}\\
& \delta_j^+ \ge 0,\quad \delta_j^- \ge 0,
&& \forall j\in J.
\label{eq:ord-asf-vcal-deviation-nonnegative}
\end{align}

The calibrated weighted reference parameters are then defined as
\begin{equation}
\widetilde v_j=\widehat v_j+\kappa,
\qquad \forall j\in J.
\label{eq:ord-asf-vcal-calibrated-v}
\end{equation}

The objective function \eqref{eq:ord-asf-vcal-obj} minimises the total
deviation between the shifted values of the parameters
\(\widehat v_j\) and the corresponding ideal-point-based target values.
The variables \(\delta_j^+\) and \(\delta_j^-\) are non-negative
deviation variables and measure, respectively, positive and negative
discrepancies. The set of constraints
\eqref{eq:ord-asf-vcal-target} defines the relation between the shifted
value \(\widehat v_j+\kappa\) and the target value
\(\widehat w_jg_j^*\). The deviation variables \(\delta_j^+\) and
\(\delta_j^-\) allow possible discrepancies between these two
quantities, with constraints
\eqref{eq:ord-asf-vcal-deviation-nonnegative} imposing the
non-negativity of the deviation variables.

\textcolor{black}{The proposed framework can also be combined with the direct elicitation of a reference point. Indeed, if the decision maker is able to specify a reference point
\[
\mathbf r^{DM}=(r^{DM}_1,\ldots,r^{DM}_m),
\]
this information can be incorporated into the ordinal regression model in different ways. The simplest possibility consists in fixing the reference point to
$\mathbf r^{DM}$ and estimating only the criteria weights. More generally,
the reference point may be treated as a decision variable while introducing
an additional objective aimed at keeping it as close as possible to the one
provided by the decision maker. For example, adopting the $L_1$ distance,
one may minimize
\[
\sum_{j\in J}\delta_j,
\]
subject to
\[
\delta_j\ge r_j-r^{DM}_j,\qquad
\delta_j\ge r^{DM}_j-r_j,
\qquad j\in J.
\]
This objective can be combined with the ordinal regression objective,
leading, for example, to the optimization problem
\[
\min
\sum_{a\in A^R}
(\sigma^+(a)+\sigma^-(a))
+
\lambda
\sum_{j\in J}\delta_j,
\]
where $\lambda\ge0$ controls the trade-off between faithfully reproducing
the decision maker's preference information and remaining close to the
reference point directly elicited from the decision maker. This extension provides a natural bridge between approaches based on the
direct specification of reference points and approaches relying on preference
disaggregation. In particular, it allows the model to exploit any reference-point
information that the decision maker is willing to provide, while automatically
adjusting it whenever this improves the consistency with the observed
preference statements.}

\subsection{ORD-ASF based Alternative Generation with Estimated Preference Parameters}

Once the estimated preference parameters have been obtained from the previous model, they are kept fixed and used in a subsequent solution-generation model. This model aims to identify a new feasible alternative that maximises the ASF-based utility function under the same estimated weights and reference values. Therefore, the parameters obtained from the preference learning phase are treated as fixed input parameters used to evaluate feasible alternatives. The model can be applied iteratively. At each iteration, the ASF-based utility function is maximized, while the alternatives previously evaluated by the DM are excluded from the feasible set.

Let $\mathcal{A}$ be the set of feasible alternatives that can be generated, and let $a \in \mathcal{A}$ denote a new candidate alternative. Moreover, let $a^{(k)}$, $k=1,\dots,\theta-1$, denote the alternatives already generated in the previous iterations.

Given the criterion weights $\widehat w_j$, $j \in J$, the corresponding weighted reference parameters $\widehat v_j$, $j \in J$, and the calibrated weighted reference parameters $\widetilde v_j$, $j \in J$, obtained from the \textbf{ORD-ASF} or the \textbf{ORD-ASF-WB} model, the alternative-generation model, referred to as \textbf{ORD-ASF-AGM}, is formulated as follows:

\begin{align}
\max_{a\in A,\,t}\quad
&
\rho\sum_{j\in J}
\left(
 \widehat w_j g_j(a)-\widetilde v_j
\right)
-
(1-\rho)t
\label{eq:asgm_obj}
\\[2mm]
\text{s.t.}\quad
&
t\geq
\widetilde v_j-\widehat w_jg_j(a),
&&
\forall j\in J,
\label{eq:asgm_t}
\\[1mm]
&
a\neq a^{(k)},
&&
k=1,\ldots,\theta-1,
\label{eq:asgm_exclusion}
\\[1mm]
&
a\in A,
\label{eq:asgm_feasible}
\\[1mm]
&
t\in\mathbb{R}.
\label{eq:asgm_tdomain}
\end{align}

The auxiliary variable $t$ represents the maximum weighted deviation
from the calibrated weighted reference vector. Indeed, since $t$ enters
the maximisation objective with the negative coefficient $-(1-\rho)$,
at optimality it takes the value
$
 t=
 \max_{j\in J}
 \left\{
  \widetilde v_j-\widehat w_jg_j(a)
 \right\}.
$ The exclusion constraints~\eqref{eq:asgm_exclusion} prevent the model
from returning alternatives already considered in previous iterations.
The model can therefore be applied iteratively. At iteration $\theta$,
the parameters $\widehat w$ and $\widetilde v$ remain fixed, while the
alternatives generated in iterations $1,\ldots,\theta-1$ are excluded
from the feasible set. The newly generated alternative is then either
presented to the DM or added to the exclusion set before solving the
next iteration.

\subsection{The ORD-ASF framework}

Let us close this section highlighting that the   models introduced can be integrated to create the so called ORD-ASF framework. The first model, \textbf{ORD-ASF}, estimates the preference parameters of the ASF-based utility function from the ordinal preference information provided by the DM. In particular, it identifies the criterion weights and the reference values 
that best reproduce the preference statements expressed over a set of reference 
alternatives.

Starting from this basic formulation, two variants are introduced to support the 
preference elicitation process. The model \textbf{ORD-ASF-WB} extends 
\textbf{ORD-ASF} by incorporating additional information on the admissible range 
of the criterion weights. This model can be used when the weights estimated by 
\textbf{ORD-ASF} appear to be poorly balanced or not fully consistent with the 
DM's view of the relative importance of the criteria. 

The model \textbf{ORD-ASF-Vcal} is instead devoted to the calibration and 
analysis of the reference values of the ASF-based utility function. It can be 
used to present the estimated reference points to the DM and to discuss their 
meaning in relation to the performances of the considered alternatives.  Therefore, \textbf{ORD-ASF-Vcal} provides an additional interaction 
mechanism between the analyst and the DM, allowing the reference points to be 
 better interpreted.

Finally, once the preference parameters have been estimated and possibly 
refined, the model \textbf{ORD-ASF-AGM} is introduced as an alternative 
generation model. This model keeps fixed the weights and reference values 
obtained from the previous preference learning phase and searches for new 
feasible alternatives that maximise the  ASF-based utility function. The 
model can be applied iteratively: at each iteration, the ASF-based utility 
function is maximised, while the alternatives previously evaluated by the DM are 
excluded from the feasible set. In this way, \textbf{ORD-ASF-AGM} makes it 
possible to identify additional alternatives that are consistent with the learned 
preference model and that may achieve a better value of the ASF-based utility 
function.

The overall framework can therefore be used in an iterative way. A possible 
procedure starts with \textbf{ORD-ASF}, which estimates the initial weights and 
reference values from the available ordinal preference information. If the 
estimated weights are not satisfactory or appear unbalanced, \textbf{ORD-ASF-WB} 
can be used to impose suitable bounds and obtain a revised set of weights. Then, 
\textbf{ORD-ASF-Vcal} can be employed to analyse and discuss the reference 
values with the DM. After this calibration phase, \textbf{ORD-ASF-AGM} can be 
run to generate new promising alternatives under the same estimated preference 
parameters. These newly generated alternatives can then be submitted to the DM 
for further evaluation, producing new preference information. The whole process 
can be repeated by returning to \textbf{ORD-ASF} and updating the estimated 
weights and reference values based on a new set of alternatives 
and preference statements. 

\begin{figure}[H]
    \centering
    \includegraphics[
        width=\textwidth,
        height=1.2\textheight,
        keepaspectratio
    ]{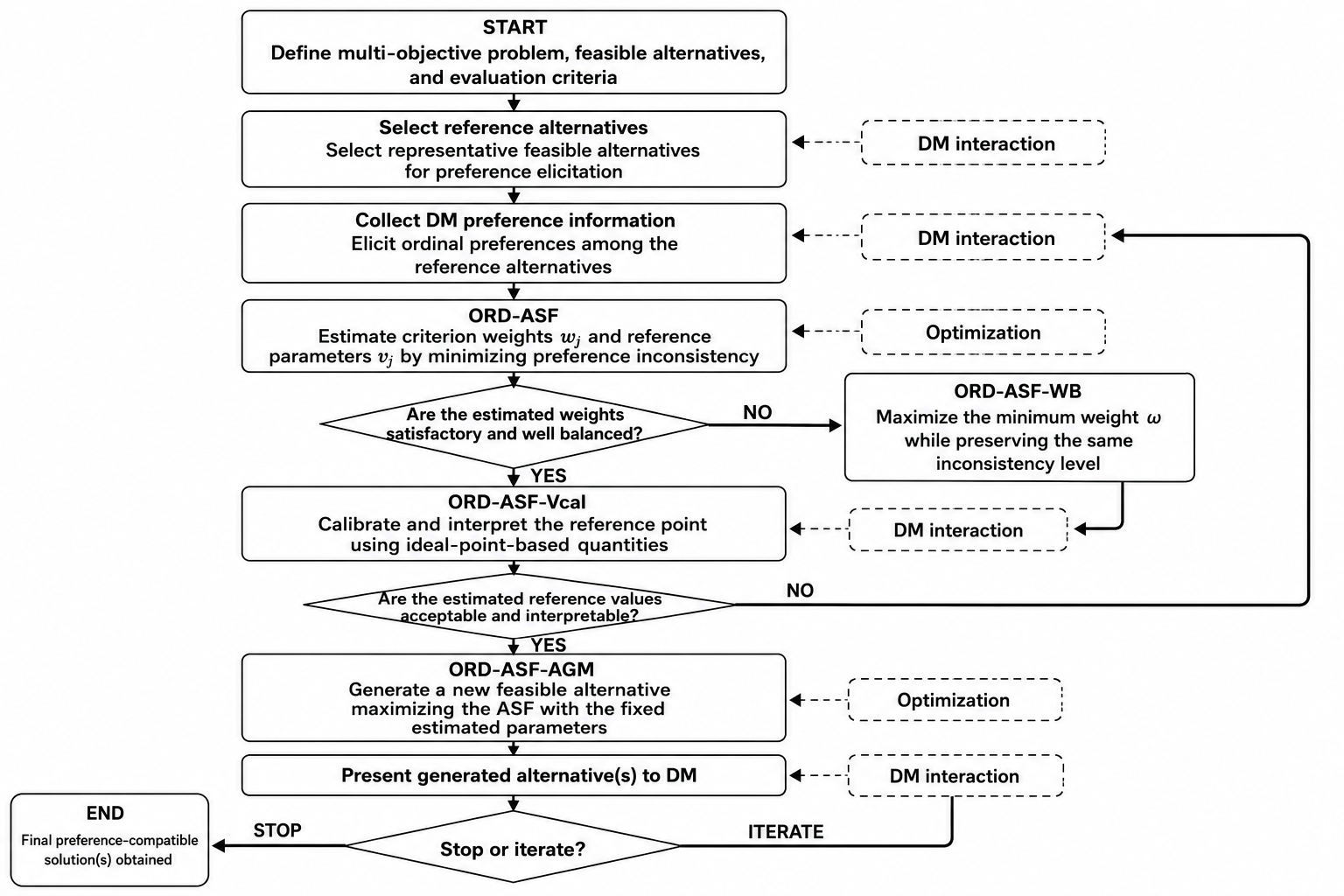}
    \caption{ORD-ASF framework}
    \label{fig:ord-asf-framework}
\end{figure}

As summarised in Figure~\ref{fig:ord-asf-framework}, the procedure starts by selecting representative alternatives and collecting the DM's ordinal preferences, which are used by ORD-ASF to estimate criterion weights and reference values. If the weights are unsatisfactory, ORD-ASF-WB improves their balance, while ORD-ASF-Vcal supports the validation and interpretation of the estimated reference values. Once the parameters are accepted, ORD-ASF-AGM generates new promising alternatives, which are evaluated by the DM and either accepted as final solutions or used to update the preference model in a new iteration.

\section{Illustrative application: the ORD-ASF framework for the multiobjective knapsack problem}\label{sec:example}

\textcolor{black}{To illustrate the proposed ordinal regression approach and its ability to induce the weights and the reference point of an ASF, we consider a multiobjective knapsack problem. The knapsack problem is adopted as an illustrative application because it provides a simple and intuitive combinatorial optimisation setting in which the alternatives can be explicitly represented and evaluated according to multiple objectives. Importantly, the choice of the knapsack problem is not restrictive for the proposed preference elicitation framework. As discussed in Section~\ref{sec:method}, the same approach can be applied to other multiobjective optimisation problems by appropriately defining the feasible set and the objective functions. The purpose of the present section is therefore not to introduce a knapsack-specific methodology, but rather to provide a concrete numerical illustration of the proposed general framework.
}
We consider the set of criteria, denoted by $J=\{1,\ldots, m\}$, and  a set of items, denoted by $I=\{1,\ldots,n\}$, that can be inserted into the knapsack. For each item $i \in I$ and each criterion $j \in J$, let $g_{ij}$ denote the contribution of item $i$ to criterion $j$, and let $c_i$ denote the cost of inserting item $i$ into the knapsack. A total budget $B$ is fixed. Each alternative $a\subseteq I$ is a feasible knapsack, that is, a subset of items from $I$ that satisfy the budget constraint, i.e. alternative $a$ is feasible if $\sum_{i \in a} c_i \leq B$. The evaluation of alternative $a$ to criterion $j$ is then given by $g_j(a)=\sum_{i \in a} g_{ij}$ where $a \subseteq I$ represents the set of items selected in the knapsack. 
Each alternative $a$ can be represented by a binary vector $x=(x_1,\ldots,x_{n})$, where $x_i=1$ if item $i$ belongs to alternative $a$, and $x_i=0$ otherwise. Thus, every feasible alternative $a$ is uniquely represented by a binary vector $x$. Therefore, the set of selected items is given by $ a=\{i \in I : x_i=1\}.$
The multiobjective knapsack problem can then be formulated as follows:
\[
\begin{aligned}
\max_{a \subseteq I} \quad 
& \bigl(g_1(a),\ldots,g_m(a)\bigr) \\
\text{s.t.}\quad 
& \sum_{i \in a} c_i \leq B.
\end{aligned}
\]

\noindent where, for each criterion $j\in J$,
\[
g_j(a)=\sum_{i\in a} g_{ij}.
\]

 To represent the DM's preferences over feasible alternatives, we use a Cobb-Douglas utility function. Let 
\[
g(a)=\bigl(g_1(a),\ldots,g_m(a)\bigr)
\]
denote the vector of criterion values associated with a feasible alternative $a$. We assume that $g_j(a)>0$ for each $j\in J=\{1,\ldots,m\}$ and that the DM evaluates each feasible alternative through the utility function
\[
U(g(a))=\prod_{j=1}^{m} g_j(a)^{\alpha_j},
\]
where $\alpha_j\ge 0$ for all $j\in J$ and $\sum_{j=1}^{m}\alpha_j=1$.

Equivalently, since the utility depends on the feasible alternative $a$ through the criterion values $g_1(a),\ldots,g_m(a)$, we can write
\[
U(a)=\prod_{j=1}^{m} g_j(a)^{\alpha_j}.
\]
The parameters $\alpha_j$, $j\in J$, express the relative importance assigned by the DM to the different criteria. The utility function is then used as a preference model to compare feasible alternatives. The evaluation $g_j(a)$ represent the evaluations that will be used in the \textbf{ORD-ASF} models introduced in the previous section. Let us note that, although we selected a Cobb–Douglas utility function for the baseline analysis, the framework is not restricted to this specific functional form. Other utility functions could also be used, depending on the assumptions one wants to impose on preferences.

It is worth presenting the formal formulation of the \textbf{ORD-ASF-AGM} model. For each item $i \in I$, let the binary variable $x_i$ indicate whether item $i$ is selected in alternative $a$. Specifically, $x_i=1$ if item $i$ belongs to knapsack $a$, and $x_i=0$ otherwise.
Hence, each feasible alternative $a \subseteq I$ is represented by the binary vector $x=(x_1,\ldots,x_n)$, and its value with respect to criterion $j \in J$ is given by
$
g_j(x)=\sum_{i\in I} g_{ij}x_i.
$
Let $\widehat v_j$, $j \in J$, denote the weighted reference parameters obtained from the  \textbf{ORD-ASF} or the \textbf{ORD-ASF-WB} model, and let $\widehat w_j$, $j \in J$, denote the corresponding estimated criterion weights. Moreover, let $x^{(k)}=(x_1^{(k)},\ldots,x_n^{(k)})$, $k=1,\ldots,\theta-1$, denote the binary vectors associated with the alternatives generated in the previous iterations.

The alternative generation model \textbf{ORD-ASF-AGM} for the multiobjective knapsack problem, can be formulated as follows:

\begin{align}
\max_{x,t} \quad
& \rho \sum_{j \in J}(
\widehat w_j \sum_{i \in I} g_{ij}x_i -\widehat v_j)
-
(1-\rho)t
\label{eq:kp_obj_simplified} \\[2mm]
\text{s.t.} \quad
& t \geq \widehat v_j - \widehat w_j \sum_{i \in I} g_{ij}x_i,
&& \forall j \in J,
\label{eq:kp_t_simplified} \\[1mm]
& \sum_{i \in I} c_i x_i \leq B,
\label{eq:kp_budget_simplified} \\[1mm]
& \sum_{i:\,x_i^{(k)}=1} \left(1-x_i\right)
+
\sum_{i:\,x_i^{(k)}=0} x_i
\geq 1,
&& k=1,\ldots,\theta-1.\label{eq:ord-asf-agm-kp-exclusion}\\
& x_i \in \{0,1\},
&& \forall i \in I,
\label{eq:kp_binary_simplified} \\[1mm]
& t \in \mathbb{R}.
\label{eq:kp_tdomain_simplified}
\end{align}

The objective function \eqref{eq:kp_obj_simplified} is composed of the two terms to a scalarising function. Constraints \eqref{eq:ord-asf-agm-kp-exclusion} prevent the generation of alternatives that have already been obtained in previous iterations. Constraints \eqref{eq:asgm_exclusion} ensure that the newly generated knapsack does not contain all the items selected in any previously generated alternative. Therefore, for each previously generated solution \(x^{(k)}\), at least one item selected in \(x^{(k)}\) must not be selected in the new solution.
Constraint \eqref{eq:kp_budget_simplified} imposes the  limit on the budget while constraints \eqref{eq:kp_binary_simplified} the binary nature of the decision variables.

\subsection {Didactic example}

Let us assume that we have an instance of the knapsack problem with 3 criteria, $J=\{1,\ldots,3\}$, 10 items $I=\{1,\ldots,10\}$, a total budget $B=10$ and the  data reported in Table~\ref{tab:items}. All criteria are maximised.

\begin{table}[htbp]
\centering
\caption{Items of the instance for the didactic example.}
\label{tab:items}
\begin{tabular}{ccccc}
\toprule
Item $i$ & $g_{i1}$ & $g_{i2}$ & $g_{i3}$ & $c_i$\\
1&5&2&5&2.3\\
2&6&6&5&2.6\\
3&2&6&1&2.1\\
4&1&5&4&1.9\\
5&4&6&1&1.8\\
6&6&4&2&1.6\\
7&1&4&2&1.9\\
8&2&6&4&2.0\\
9&4&1&3&1.3\\
10&3&3&5&1.6\\
\bottomrule
\end{tabular}
\end{table}

Enumerating the $2^{10}$ binary vectors gives $584$ feasible alternatives. The
non-dominated set contains $14$ alternatives; the best one according to
Cobb-Douglas utility Function is
\begin{equation*}
 x^\star=\vect{1,1,0,0,1,1,0,0,0,1},
 \qquad g(x^\star)=\vect{24,21,18},
 \qquad U(x^\star)=21.7684.
\end{equation*}
Let us note that this solution is used only as an external benchmark to verify the generated
alternatives; it is not supplied to the ordinal-regression model. In the didactic example, the criteria are not normalised because their coefficients
are generated within the same numerical range from the same uniform distribution as conventional in randomly generated instances for multiobjective knapsack problems
\citep[e.g.,][]{MesquitaCunha2023}. Therefore, the criteria are comparable by construction.
In applications
involving criteria expressed in different units or characterised by substantially different
ranges, the objective values should be suitably normalised before applying the proposed
ordinal-regression framework.

We then select a subset of four non-dominated knapsacks to form the set $A^R$ of representative alternatives. These alternatives are used to express the DM's preference information. In this example, the ordering of the representative alternatives is obtained according to the following Cobb-Douglas utility function \[ U(a)=g_1(a)^{0.5}g_2(a)^{0.3}g_3(a)^{0.2}. \] and we report them in Table~\ref{tab:selected-reference-alternatives}, sorted in decreasing order of Cobb-Douglas utility and, therefore, according to the DM's preferences. For compactness, each alternative $a$ is reported through the corresponding selected-item set $a=\{i \in I : x_i=1\}.$ 

\begin{table}[htbp]
\centering
\caption{Selected reference alternatives.}
\label{tab:selected-reference-alternatives}
\begin{tabular}{c c c c c c}
\toprule
Reference alternative & Selected items & $g_1(a^{(r)})$ & $g_2(a^{(r)})$ & $g_3(a^{(r)})$ & $U(a^{(r)})$ \\
\midrule
$a^{(1)}$ & $\{1,2,6,8,9\}$      & 23 & 19 & 19 & 20.9045 \\
$a^{(2)}$ & $\{2,4,5,6,8\}$      & 19 & 27 & 16 & 20.3990 \\
$a^{(3)}$ & $\{2,4,6,8,10\}$     & 18 & 24 & 20 & 20.0404 \\
$a^{(4)}$ & $\{1,2,8,9,10\}$     & 20 & 18 & 22 & 19.7506 \\
\bottomrule
\end{tabular}
\end{table}

Thus, the preference information provided by the DM is expressed by the following complete order over the four reference alternatives: $a^{(1)} \succ a^{(2)} \succ a^{(3)} \succ a^{(4)}.$

We also compute the ideal objective vector by solving three single-objective knapsack problems, one for each criterion. More precisely, for each $j \in J=\{1,2,3\}$, we compute
$g^*_j=\max_{a \in A} g_j(a),$ where $A$ denotes the set of feasible alternatives satisfying the budget constraint, namely $A=\left\{a\subseteq I:\sum_{i\in a} c_i \leq B\right\}.$
For this instance, the resulting ideal objective vector is $ g^*=(g_1^*, g_2^*, g_3^*)=(25,27,22).$

Solving the \textbf{ORD-ASF} model with an optimisation solver, based on
the preference information provided over the set $A^R$,  imposing the
strict preference margin $\varepsilon=0.1$, yields an optimal
inconsistency level equal to zero. The parameter $\varepsilon>0$ represents the minimum required separation
between the utilities of two consecutively ranked alternatives. Its
introduction prevents alternatives declared strictly preferred by the
DM from receiving the same utility value.
\textcolor{black}{Let us note that in this particular instance, the same optimal inconsistency level
is obtained for all values of \(\rho\)  considered in the grid, namely
\(\rho \in \{0.05,0.10,\ldots,1.00\}\); All the considered values of $\rho$ yield zero preference inconsistency. However, the estimated parameters and the alternatives generated in the subsequent phase may still depend on $\rho$. We set $\rho=0.3$ for illustrative purposes in
equation~\eqref{eq:ord-asf-value}.}

One optimal solution returned by the solver is
$ w=(w_1,w_2,w_3)
 = \vect{0.2429,\,0.4333,\,0.3238},
$
with weighted reference parameters
$
 v=(v_1,v_2,v_3)
 =
 \vect{0,\,0.4333,\,0.5667}.
$

The utilities obtained from the \textbf{ORD-ASF} model
are reported in Table~\ref{tab:initial_regression}. They reproduce the
strict preference order provided by the DM with the required margin
$\varepsilon=0.1$. Since the optimal solution of the \textbf{ORD-ASF} model may
 not be unique, the reported weight vector should be interpreted as
one model compatible with the preference information, rather than as
the unique representation of the DM's preferences.

\begin{table}[htbp]
\centering
\caption{Utilities of the reference alternatives under an initial
ORD-ASF optimal solution.}
\label{tab:initial_regression}
\begin{tabular}{cccccc}
\toprule
Alternative &  $g_1$ & $g_2$ & $g_3$ & $U(a)$\\
\midrule
$a^{(1)}$   & 23 & 19 & 19 & 9.9014\\
$a^{(2)}$  & 19 & 27 & 16 & 9.6786\\
$a^{(3)}$  & 18 & 24 & 20 & 9.4343\\
$a^{(4)}$  & 20 & 18 & 22 & 9.3343\\
\bottomrule
\end{tabular}
\end{table}
\textcolor{black}{
Since the optimal solution of ORD-ASF is not necessarily unique, the analyst may use
a secondary criterion to select a representative parameter vector among the solutions
that preserve the same inconsistency level. In particular, instead of fixing the strict
preference margin \(\epsilon\) a priori, one may maximise the minimum separation between
the corrected utilities of two consecutive reference alternatives, while keeping the
optimal inconsistency level obtained by ORD-ASF unchanged, as described in the second extension of the model.
Let \(\sigma^*=0\) be the
optimal inconsistency level obtained in the previous ORD-ASF step. We then search,
among the solutions satisfying
$
\sum_{r=1}^{|A^R|}
\left(\sigma^+(a^{(r)})+\sigma^-(a^{(r)})\right)
\leq \sigma^*,$
for the one that maximises the minimum value of
$
 U(a^{(r)})- U(a^{(r+1)}),
\qquad r=1,\ldots,|A^R|-1.
$
For the same value \(\rho=0.3\), the maximum-separation refinement gives a minimum
corrected-utility distance equal to
\[
\epsilon^*=0.2265,
\]
while the total deviation remains equal to zero. One corresponding set of parameters is
\[
w=[0.7454,\,0.1252,\,0.1294],
\qquad
v=[0,\,1,\,0],
\]
which gives the reference vector
\[
r=[0,\,7.9852,\,0].
\]
The resulting utilities of the four reference alternatives are
\[
U(a^{(1)})=7.5600,\qquad
U(a^{(2)})=7.3335,\qquad
U(a^{(3)})=7.1070,\qquad
U(a^{(4)})=6.8804.
\]
}

\textcolor{black}{
We also consider the extension in which the Tchebycheff and weighted-sum components are
allowed to have two distinct weight vectors. 
In this formulation, the parameter \(\rho\) is not fixed a priori but is recovered from
the scale of the weighted-sum component. Solving this variant on the same preference
information gives an optimal inconsistency level equal to zero. The resulting Tchebycheff
weights are
\[
w^T=[0.3231,\,0.3385,\,0.3385],
\]
while the weighted-sum weights are
\[
\widetilde w^S=[0,\,0.1244,\,0].
\]
Therefore,
\[
\mu=\sum_{j\in J}\widetilde w^S_j = 0.1244,
\]
and the corresponding value of \(\rho\) is
\[
\rho=\frac{\mu}{1+\mu}
     =\frac{0.1244}{1+0.1244}
     =0.1106.
\]
The normalised weighted-sum vector is
\[
w^S=\frac{\widetilde w^S}{\mu}=[0,\,1,\,0].
\]
The weighted reference vector is
\[
v^T=[1,\,0,\,0],
\]
which gives the reference point
\[
r=\left[\frac{1}{0.3231},\,0,\,0\right]=[3.0952,\,0,\,0].
\]
The utilities of the four reference alternatives are
\[
U(a^{(1)})=8.7936,\quad
U(a^{(2)})=8.4962,\quad
U(a^{(3)})=7.8000,\quad
U(a^{(4)})=7.7000.
\]
Thus, the strict preference order is reproduced without deviations, since
\[
\sigma^+=\sigma^-=[0,\,0,\,0,\,0].
\]
This result shows that allowing separate weights for the two ASF components provides
an alternative representation of the same preference information. In particular, the
Tchebycheff component uses a balanced weight vector, whereas the compensatory
weighted-sum component is concentrated on the second criterion, with an endogenous
relative contribution \(\rho=0.1106\).
}

Although the initial solution assigns a positive weight to every
criterion, the \textbf{ORD-ASF-WB} model is subsequently applied to
select, among all solutions preserving the optimal inconsistency level
$\sigma^*=0$, one in which the minimum criterion weight is as large
as possible. The resulting balanced weights are
\begin{equation}
 \widehat w
 =
 \vect{0.3278,\,0.3278,\,0.3444},
 \qquad
 \omega^*=0.3278.
 \label{eq:balanced_weights}
\end{equation}
A corresponding optimal weighted reference vector is
\begin{equation}
 \widehat v
 =
 \vect{0.5620,\,0.4380,\,0}.
 \label{eq:balanced_v}
\end{equation}

Table~\ref{tab:wb_check} verifies that the strict preference order defined on $A^R$ remains compatible with the estimated parameters.

\begin{table}[htbp]
\centering
\caption{Reference alternatives under the balanced-weight solution.}
\label{tab:wb_check}
\begin{tabular}{ccccccc}
\toprule
Alternative  & $g_1$ & $g_2$ & $g_3$ &  $U(a)$\\
\midrule
$a^{(1)}$ 
 & 23 & 19 & 19 &  10.1467\\
$a^{(2)}$
 & 19 & 27 & 16 &  10.0338\\
$a^{(3)}$
 & 18 & 24 & 20 &  9.9338\\
$a^{(4)}$ 
 & 20 & 18 & 22 &  9.8338\\
\bottomrule
\end{tabular}
\end{table}

The calibration model \textbf{ORD-ASF-Vcal} is solved after
the weight-balancing phase. In this way, the reference vector used in
the subsequent alternative-generation phase is calibrated consistently
with the balanced weight vector $\widehat w$.

Using the ideal objective vector
$ g^*=\vect{25,\,27,\,22},$
the model shifts the weighted reference vector $\widehat v$ by a common
scalar $\kappa$ and compares it with the ideal-point-based quantities
$\widehat w_jg_j^*$. The optimal common shift is
$
 \kappa=7.6334.$
The calibrated weighted reference vector is therefore
\begin{equation}
 \widetilde v_j
 =
 \widehat v_j+\kappa,
 \qquad j\in J,
\end{equation}
that is,
\begin{equation}
 \widetilde v
 =
 \vect{8.1954,\,8.0714,\,7.6334}.
 \label{eq:vcal_corrected_v}
\end{equation}

The corresponding deviation variables are
\begin{equation}
 \delta^+
 =
 \vect{0,\,0,\,0.0572},
 \qquad
 \delta^-
 =
 \vect{0,\,0.7796,\,0},
\end{equation}
and the optimal calibration deviation is
\begin{equation}
 \sum_{j\in J}
 \left(
 \delta_j^+ + \delta_j^-
 \right)
 =
 0.8368.
\end{equation}

The calibrated reference points are obtained as
\begin{equation}
 \widetilde r_j
 =
 \frac{\widetilde v_j}{\widehat w_j},
 \qquad j\in J,
\end{equation}
giving
\begin{equation}
 \widetilde r
 =
 \vect{25,\,24.6219,\,22.1662}.
 \label{eq:vcal_corrected_reference}
\end{equation}

This calibrated reference point should not be interpreted as the unique
behavioural representation of the DM's aspirations. Rather, it is one
representative reference point obtained from the family of parameter
vectors compatible with the stated preference information and with the
calibration criterion.

Once the weights and the reference vector have been estimated and
calibrated, they are kept fixed in the alternative-generation model
\textbf{ORD-ASF-AGM}. In particular, the generation phase uses
\[
 \rho=0.3,
 \qquad
 \widehat w
 =
 \vect{0.3278,\,0.3278,\,0.3444},
\]
and
\[
 \widetilde v
 =
 \vect{8.1954,\,8.0714,\,7.6334}.
\]

After excluding the four reference alternatives reported in
Table~\ref{tab:selected-reference-alternatives}, the first execution of
\textbf{ORD-ASF-AGM} returns
\begin{equation}
 x^{(1)}
 =
 \vect{1,1,0,0,1,1,0,0,0,1}.
\end{equation}
The corresponding objective vector is
\begin{equation}
 g(x^{(1)})
 =
 \vect{24,\,21,\,18}.
\end{equation}
For this solution,
\[
 t(x^{(1)})=1.4347,
\]
and the \textbf{ORD-ASF-AGM} objective value is
$
 -1.8893.
$

The generated alternative coincides with the external
Cobb-Douglas benchmark solution $x^\star$, whose utility is
$
 U(x^\star)=21.7684.
$
It is important to recall that the benchmark solution and its
Cobb-Douglas utility are used only for ex-post verification and are not
provided to any of the ordinal-regression or alternative-generation
models.

Since the generated alternative coincides with $x^\star$, the procedure
may stop immediately if the adopted stopping rule is to terminate when
the benchmark is recovered. Further iterations would be required 
if the analyst wished to generate and present additional
preference-compatible alternatives to the DM and the iteration process may continue solving again the \textbf{ORD-ASF} model to estimate new weights and reference points. 

Let us note that, in this example, we assume a set of parameters, such as $\rho$, $\epsilon$, and the number of solutions to be presented to the DM. These parameters may vary and may help identify the DM's different preferences.

\section{Conclusions}\label{sec:conclusions}
In this paper, we proposed a preference elicitation procedure to infer the parameters of an achievement scalarising function in multiobjective optimisation. In particular, the approach allows us to estimate both the weights associated with the objective functions and the reference point representing the aspiration levels of the DM. Instead of requiring the DM to specify these parameters directly, the method relies on holistic preference information expressed on a set of alternatives and uses ordinal regression to identify a scalarising function compatible with such preferences.

The proposed framework contributes to bridging preference learning and multiobjective optimisation. On the one hand, it reduces the cognitive effort required from the DM, who is not asked to provide precise numerical values for weights or aspiration levels. On the other hand, the inferred scalarising function can be directly embedded into an optimisation procedure, guiding the search towards regions of the Pareto front that are coherent with the expressed preferences. The application to the multiobjective knapsack problem illustrates the practical implementation of the method and shows how preference information can be translated into a scalar optimisation model. 

Several directions for future research can be considered. First, the approach should be tested on further classical multiobjective optimisation problems, such as assignment, location and portfolio selection problems. \textcolor{black}{Indeed, although the proposed approach has been illustrated on a multiobjective knapsack problem, the underlying ordinal regression framework is general and can be applied to other multiobjective optimisation settings.} Second, future work may investigate the extension of the ordinal regression procedure to other families of utility or scalarising functions, beyond achievement scalarising functions, to assess whether different preference structures can be captured more effectively. Third, the use of multiple or progressively revised weights and reference points could be explored, especially within interactive procedures in which the DM refines preferences during the optimisation process. Also, some visual tools could be implemented to help the DM visualise the evolution of the different solutions. 

Some limitations should also be acknowledged. The quality of the inferred model depends on the consistency of the preference information provided by the DM. Also, iterative procedures may require careful design to avoid excessive cognitive burden and to ensure that the interaction remains meaningful for the DM. Addressing these issues represents an important step for further developing preference-based scalarisation methods in multiobjective optimisation.

\bibliography{Sample.bib}

@article{MesquitaCunha2023,
  author  = {Mesquita-Cunha, Mariana and Figueira, Jos{\'e} Rui and Barbosa-P{\'o}voa, Ana Paula},
  title   = {New {$\varepsilon$}-constraint methods for multi-objective integer linear programming: A {Pareto} front representation approach},
  journal = {European Journal of Operational Research},
  volume  = {306},
  number  = {1},
  pages   = {286--307},
  year    = {2023},
  doi     = {10.1016/j.ejor.2022.07.044}
}

@incollection{wierzbicki1980use,
  author    = {Wierzbicki, Andrzej P.},
  title     = {The Use of Reference Objectives in Multiobjective Optimization},
  booktitle = {Multiple Criteria Decision Making Theory and Application},
  editor    = {Fandel, G. and Gal, T.},
  series    = {Lecture Notes in Economics and Mathematical Systems},
  volume    = {177},
  pages     = {468--486},
  publisher = {Springer},
  address   = {Berlin, Heidelberg},
  year      = {1980},
  doi       = {10.1007/978-3-642-48782-8_32},
}

@article{greco2008ordinal,
  author  = {Greco, Salvatore and Mousseau, Vincent and S{\l}owi{\'n}ski, Roman},
  title   = {Ordinal Regression Revisited: Multiple Criteria Ranking Using a Set of Additive Value Functions},
  journal = {European Journal of Operational Research},
  volume  = {191},
  number  = {2},
  pages   = {416--436},
  year    = {2008},
  doi     = {10.1016/j.ejor.2007.08.013},
}

@book{miettinen1999nonlinear,
  author    = {Miettinen, Kaisa},
  title     = {Nonlinear Multiobjective Optimization},
  series    = {International Series in Operations Research \& Management Science},
  volume    = {12},
  publisher = {Springer},
  address   = {Boston, MA},
  year      = {1999},
  doi       = {10.1007/978-1-4615-5563-6},
}

@article{miettinen2002scalarising,
  author  = {Miettinen, Kaisa and M{\"a}kel{\"a}, Marko M.},
  title   = {On Scalarizing Functions in Multiobjective Optimization},
  journal = {OR Spectrum},
  volume  = {24},
  number  = {2},
  pages   = {193--213},
  year    = {2002},
  doi     = {10.1007/s00291-001-0092-9},
}

@article{nikulin2012new,
  author  = {Nikulin, Yury and Miettinen, Kaisa and M{\"a}kel{\"a}, Marko M.},
  title   = {A New Achievement Scalarizing Function Based on Parameterization in Multiobjective Optimization},
  journal = {OR Spectrum},
  volume  = {34},
  number  = {1},
  pages   = {69--87},
  year    = {2012},
  doi     = {10.1007/s00291-010-0224-1},
}

@article{luque2009incorporating,
  author  = {Luque, Mariano and Miettinen, Kaisa and Eskelinen, Petri and Ruiz, Francisco},
  title   = {Incorporating Preference Information in Interactive Reference Point Methods for Multiobjective Optimization},
  journal = {Omega},
  volume  = {37},
  number  = {2},
  pages   = {450--462},
  year    = {2009},
  doi     = {10.1016/j.omega.2007.06.001},
}

@article{klamroth2008integrating,
  author  = {Klamroth, Kathrin and Miettinen, Kaisa},
  title   = {Integrating Approximation and Interactive Decision Making in Multicriteria Optimization},
  journal = {Operations Research},
  volume  = {56},
  number  = {1},
  pages   = {222--234},
  year    = {2008},
  doi     = {10.1287/opre.1070.0425},
}

@incollection{greco2010robust,
  author    = {Greco, Salvatore and S{\l}owi{\'n}ski, Roman and Figueira, Jos{\'e} Rui and Mousseau, Vincent},
  title     = {Robust Ordinal Regression},
  booktitle = {Trends in Multiple Criteria Decision Analysis},
  editor    = {Ehrgott, Matthias and Figueira, Jos{\'e} Rui and Greco, Salvatore},
  series    = {International Series in Operations Research \& Management Science},
  volume    = {142},
  pages     = {241--283},
  publisher = {Springer},
  address   = {Boston, MA},
  year      = {2010},
  doi       = {10.1007/978-1-4419-5904-1_9},
}

@article{corrente2013robust,
  author  = {Corrente, Salvatore and Greco, Salvatore and Kadzi{\'n}ski, Mi{\l}osz and S{\l}owi{\'n}ski, Roman},
  title   = {Robust Ordinal Regression in Preference Learning and Ranking},
  journal = {Machine Learning},
  volume  = {93},
  number  = {2--3},
  pages   = {381--422},
  year    = {2013},
  doi     = {10.1007/s10994-013-5365-4},
}

@article{greco2014robust,
  author  = {Greco, Salvatore and Mousseau, Vincent and S{\l}owi{\'n}ski, Roman},
  title   = {Robust Ordinal Regression for Value Functions Handling Interacting Criteria},
  journal = {European Journal of Operational Research},
  volume  = {239},
  number  = {3},
  pages   = {711--730},
  year    = {2014},
  doi     = {10.1016/j.ejor.2014.05.022},
}

@article{barbati2024deck,
  author  = {Barbati, Maria and Greco, Salvatore and Lami, Isabella M.},
  title   = {The Deck-of-cards-based Ordinal Regression Method and Its Application for the Development of an Ecovillage},
  journal = {European Journal of Operational Research},
  volume  = {319},
  number  = {3},
  pages   = {845--861},
  year    = {2024},
  doi     = {10.1016/j.ejor.2024.07.010},
}

@article{greco2025fifty,
  author  = {Greco, Salvatore and S{\l}owi{\'n}ski, Roman and Wallenius, Jyrki},
  title   = {Fifty Years of Multiple Criteria Decision Analysis: From Classical Methods to Robust Ordinal Regression},
  journal = {European Journal of Operational Research},
  volume  = {323},
  number  = {2},
  pages   = {351--377},
  year    = {2025},
  doi     = {10.1016/j.ejor.2024.07.038},
}

@article{wierzbicki1982mathematical,
  author  = {Wierzbicki, Andrzej P.},
  title   = {A Mathematical Basis for Satisficing Decision Making},
  journal = {Mathematical Modelling},
  volume  = {3},
  number  = {5},
  pages   = {391--405},
  year    = {1982},
  doi     = {10.1016/0270-0255(82)90038-0},
}

@article{wierzbicki1986completeness,
  author  = {Wierzbicki, Andrzej P.},
  title   = {On the Completeness and Constructiveness of Parametric Characterizations to Vector Optimization Problems},
  journal = {OR Spectrum},
  volume  = {8},
  number  = {2},
  pages   = {73--87},
  year    = {1986},
  doi     = {10.1007/BF01719738},
}

@techreport{kallio1980implementation,
  author      = {Kallio, Markku and Lewandowski, Andrzej and Orchard-Hays, William},
  title       = {An Implementation of the Reference Point Approach for Multiobjective Optimization},
  number      = {WP-80-035},
  address     = {Laxenburg, Austria},
  institution = {International Institute for Applied Systems Analysis},
  year        = {1980},
}

@article{miettinen2006synchronous,
  author  = {Miettinen, Kaisa and M{\"a}kel{\"a}, Marko M.},
  title   = {Synchronous Approach in Interactive Multiobjective Optimization},
  journal = {European Journal of Operational Research},
  volume  = {170},
  number  = {3},
  pages   = {909--922},
  year    = {2006},
  doi     = {10.1016/j.ejor.2004.07.052},
}

@article{miettinen2006experiments,
  author  = {Miettinen, Kaisa and M{\"a}kel{\"a}, Marko M. and Kaario, Katja},
  title   = {Experiments with Classification-Based Scalarizing Functions in Interactive Multiobjective Optimization},
  journal = {European Journal of Operational Research},
  volume  = {175},
  number  = {2},
  pages   = {931--947},
  year    = {2006},
  doi     = {10.1016/j.ejor.2005.06.019},
}

@inproceedings{wierzbicki2007reference,
  author    = {Wierzbicki, Andrzej P.},
  title     = {Reference Point Approaches and Objective Ranking},
  booktitle = {Practical Approaches to Multi-Objective Optimization},
  series    = {Dagstuhl Seminar Proceedings (DagSemProc)},
  volume    = {6501},
  pages     = {1--20},
  publisher = {Schloss Dagstuhl -- Leibniz-Zentrum f{\"u}r Informatik},
  address   = {Dagstuhl, Germany},
  year      = {2007},
  doi       = {10.4230/DagSemProc.06501.2},
  url       = {https://drops.dagstuhl.de/entities/document/10.4230/DagSemProc.06501.2},
}

@incollection{miettinen2008introduction,
  author    = {Miettinen, Kaisa and Ruiz, Francisco and Wierzbicki, Andrzej P.},
  title     = {Introduction to Multiobjective Optimization: Interactive Approaches},
  booktitle = {Multiobjective Optimization: Interactive and Evolutionary Approaches},
  series    = {Lecture Notes in Computer Science},
  volume    = {5252},
  pages     = {27--57},
  publisher = {Springer},
  address   = {Berlin, Heidelberg},
  year      = {2008},
  doi       = {10.1007/978-3-540-88908-3_2},
}

@article{ruiz2008additive,
  author  = {Ruiz, Francisco and Luque, Mariano and Miguel, Francisco and Cabello, Jos{'e} M.},
  title   = {An Additive Achievement Scalarizing Function for Multiobjective Programming Problems},
  journal = {European Journal of Operational Research},
  volume  = {188},
  number  = {3},
  pages   = {683--694},
  year    = {2008},
}

@article{figueira2010parallel,
  author  = {Figueira, Jos{\'e} Rui and Liefooghe, Arnaud and Talbi, El-Ghazali and Wierzbicki, Andrzej P.},
  title   = {A Parallel Multiple Reference Point Approach for Multi-Objective Optimization},
  journal = {European Journal of Operational Research},
  volume  = {205},
  number  = {2},
  pages   = {390--400},
  year    = {2010},
  doi     = {10.1016/j.ejor.2009.12.027},
}

@article{luque2012two,
  author  = {Luque, Mariano and Miettinen, Kaisa and Ruiz, Ana B. and Ruiz, Francisco},
  title   = {A Two-Slope Achievement Scalarizing Function for Interactive Multiobjective Optimization},
  journal = {Computers \& Operations Research},
  volume  = {39},
  number  = {7},
  pages   = {1673--1681},
  year    = {2012},
  doi     = {10.1016/j.cor.2011.10.002},
}

@article{ruiz2015preference,
  author  = {Ruiz, Ana B. and Saborido, Rafael and Luque, Mariano},
  title   = {A Preference-Based Evolutionary Algorithm for Multiobjective Optimization: The Weighting Achievement Scalarizing Function Genetic Algorithm},
  journal = {Journal of Global Optimization},
  volume  = {62},
  number  = {1},
  pages   = {101--129},
  year    = {2015},
  doi     = {10.1007/s10898-014-0214-y},
}

@inproceedings{ishibuchi2010indicator,
  author    = {Ishibuchi, Hisao and Tsukamoto, Noritaka and Sakane, Yusuke and Nojima, Yusuke},
  title     = {Indicator-Based Evolutionary Algorithm with Hypervolume Approximation by Achievement Scalarizing Functions},
  booktitle = {Proceedings of the Genetic and Evolutionary Computation Conference},
  pages     = {527--534},
  publisher = {ACM},
  year      = {2010},
  doi       = {10.1145/1830483.1830578},
}

@article{aliano2021exact,
  author  = {Aliano Filho, Alfredo and Moretti, Antonio Carlos and Vaz Pato, Maria and Oliveira, Welington A.},
  title   = {An Exact Scalarization Method with Multiple Reference Points for Bi-Objective Integer Linear Optimization Problems},
  journal = {Annals of Operations Research},
  volume  = {296},
  number  = {1},
  pages   = {35--69},
  year    = {2021},
  doi     = {10.1007/s10479-019-03317-9},
}

@article{luque2015equivalent,
  author  = {Luque, Mariano and L{\'o}pez-Agudo, Luis A. and Marcenaro-Guti{\'e}rrez, Oscar D.},
  title   = {Equivalent Reference Points in Multiobjective Programming},
  journal = {Expert Systems with Applications},
  volume  = {42},
  number  = {4},
  pages   = {2205--2212},
  year    = {2015},
  doi     = {10.1016/j.eswa.2014.10.028},
}

@article{DECASTRO2026107469,
title = {On bivariate achievement scalarizing functions},
journal = {Operations Research Letters},
volume = {68},
pages = {107469},
year = {2026},
issn = {0167-6377},
doi = {https://doi.org/10.1016/j.orl.2026.107469},
author = {Philip J. {de Castro}},
}

@article{tanabe2024quality,
  author  = {Tanabe, Ryoji and Li, Ke},
  title   = {Quality Indicators for Preference-Based Evolutionary Multiobjective Optimization Using a Reference Point: A Review and Analysis},
  journal = {IEEE Transactions on Evolutionary Computation},
  volume  = {28},
  number  = {6},
  pages   = {1575--1589},
  year    = {2024},
  doi     = {10.1109/TEVC.2023.3319009},
}

@article{phelps2003interactive,
  author  = {Phelps, Steve and K{\"o}ksalan, Murat},
  title   = {An Interactive Evolutionary Metaheuristic for Multiobjective Combinatorial Optimization},
  journal = {Management Science},
  volume  = {49},
  number  = {12},
  pages   = {1726--1738},
  year    = {2003},
  doi     = {10.1287/mnsc.49.12.1726.25117},
}

@article{lokman2016interactive,
  author  = {Lokman, B. and K{\"o}ksalan, Murat and Korhonen, Pekka and Wallenius, Jyrki},
  title   = {An Interactive Algorithm to Find the Most Preferred Solution of Multi-Objective Integer Programs},
  journal = {Annals of Operations Research},
  volume  = {245},
  number  = {1--2},
  pages   = {67--95},
  year    = {2016},
  doi     = {10.1007/s10479-014-1545-2},
}

@article{karakaya2018interactive,
  author  = {Karakaya, G{"u}lsah and K{\"o}ksalan, Murat and Ahipa{\c{s}}ao{\u{g}}lu, Suat D.},
  title   = {Interactive Algorithms for a Broad Underlying Family of Preference Functions},
  journal = {European Journal of Operational Research},
  volume  = {265},
  number  = {1},
  pages   = {248--262},
  year    = {2018},
  doi     = {10.1016/j.ejor.2017.07.028},
}

@article{karakaya2021evaluating,
  author  = {Karakaya, G{"u}lsah and K{\"o}ksalan, Murat},
  title   = {Evaluating Solutions and Solution Sets under Multiple Objectives},
  journal = {European Journal of Operational Research},
  volume  = {294},
  number  = {1},
  pages   = {16--28},
  year    = {2021},
  doi     = {10.1016/j.ejor.2021.01.021},
}

@article{karakaya2023finding,
  author    = {Karakaya, G{\"u}lsah and K{\"o}ksalan, Murat},
  title     = {Finding Preferred Solutions under Weighted Tchebycheff Preference Functions for Multi-Objective Integer Programs},
  journal   = {European Journal of Operational Research},
  volume    = {308},
  number    = {1},
  pages     = {215--228},
  publisher = {Elsevier},
  year      = {2023},
  doi       = {10.1016/j.ejor.2022.11.043},
}

@article{karakaya2025improved,
  author    = {Karakaya, G{\"u}l{\c{s}}ah and K{\"o}ksalan, M},
  title     = {An improved algorithm exploiting the characteristics of a distance-based preference function to converge to preferred solutions},
  journal   = {European Journal of Operational Research},
  volume    = {330},
  number    = {2},
  pages     = {595--607},
  publisher = {Elsevier},
  year      = {2026},
  doi       = {10.1016/j.ejor.2025.08.036},
}

@article{jacquet1982assessing,
  author  = {Jacquet-Lagr{\`e}ze, {\`E}ric and Siskos, Yannis},
  title   = {Assessing a Set of Additive Utility Functions for Multicriteria Decision-Making, the {UTA} Method},
  journal = {European Journal of Operational Research},
  volume  = {10},
  number  = {2},
  pages   = {151--164},
  year    = {1982},
  doi     = {10.1016/0377-2217(82)90155-2},
}
\bibliographystyle{apalike} 
\end{document}